\documentclass[12pt]{article}         
\usepackage{amsfonts,amssymb} 
\makeatletter
\def\bdraft{\pagestyle{myheadings} 
           \textheight=10.5truein \textwidth=7.5truein \parindent=8pt
           \voffset=-1truein \topmargin=30pt \headheight=10pt \headsep=3pt
           \ifcase \@ptsize \hoffset=-1.5truein \or \hoffset=-1.35truein
                        \or \hoffset=-1.15truein \fi}
\def\quality{\textheight=240mm \textwidth=160mm \topmargin=0Truein
             \ifcase \@ptsize \hoffset=-23mm
                     \or \hoffset=-20mm \or \hoffset=-15mm \fi}
\makeatother \quality 

\def\beq#1#2{\begin{equation} \label{#1} #2 \end{equation}}
\def\bea#1{\begin{eqnarray*} #1 \end{eqnarray*}} \def\a{\!\!\!&\!\!\!\!&}
\def\n{\noindent}   \def\map{T}   
    
\def\IR{{\mathbb{R}}}  \def\IZ{{\mathbb{Z}}}  
\def\toas#1{\stackrel{#1}{\longrightarrow}}
  \def\phi{\varphi}   
  \def\mod1{\,({\rm mod\ } 1)\,}
\def\t{\tilde} 

\def\function#1{\left\{\!\!\!\begin{array}{ll} #1 \end{array} \right.}

\def\proof{\smallskip \noindent {\bf Proof. \ }}       
\def\blanksquare{\,\,\,$\sqcup\!\!\!\!\sqcap$}         
\def\qed{\hfill\blanksquare\linebreak\smallskip\par}   

\def\thname{Theorem}     \def\lmname{Lemma}      \def\prname{Proposition}
\def\dfname{Definition}  \def\crname{Corollary}  \def\rmname{Remark}

\newtheorem{theorem}{\thname}
\newtheorem{lemma}{\lmname}
\newtheorem{corollary}[lemma]{\crname}   

\newtheorem{dftn}{\dfname}[section]
\newtheorem{rmrk}[lemma]{\rmname}
\newenvironment{remark}{\begin{rmrk}\rm}{\end{rmrk}}     

             \catcode`@=11 \@addtoreset{equation}{section} \catcode`@=12

\makeatletter\def\fps@figure{htbp}\makeatother 
\newcommand\mlbscale{1pt} 
\newif\iffigs\figstrue 
\def\bline(#1,#2)(#3,#4)(#5){\put(#1,#2){\line(#3,#4){#5}}}  

\def\bfig(#1,#2)#3#4{\begin{figure} \begin{center}
    \framebox{\setlength{\unitlength}{\mlbscale}
       \iffigs \begin{picture}(#1,#2) #3 \end{picture}
       \else \begin{picture}(60,10)(0,0)
                   \put(0,0){\framebox(60,10){Figure}} \end{picture} \fi}
    \end{center} \caption{#4} \end{figure}}

\def\Bfig(#1,#2)#3#4{\begin{figure} \begin{center}
    \setlength{\unitlength}{\mlbscale}
       \iffigs \begin{picture}(#1,#2) #3 \end{picture}
       \else \begin{picture}(60,10)(0,0)
                   \put(0,0){\framebox(60,10){Figure}} \end{picture} \fi
    \end{center} \caption{#4} \end{figure}}

\def\bpic(#1,#2)#3{\setlength{\unitlength}{\mlbscale}
    \begin{picture}(#1,#2) #3 \end{picture}}

\def\*#1{#1^*}    \def\0#1{\breve#1}  \def\2#1{\acute#1}
\def\G{\Delta}    \def\b#1{\overline{#1}}
\def\intp#1{\left\lfloor#1\right\rfloor}

\def\?#1{} 

\def\negvskip{\vskip-0.5truecm}

\begin{document}
\title{Exclusion type spatially heterogeneous processes \\ in continuum}
\author{Michael Blank\thanks{
        Russian Academy of Sci., Inst. for
        Information Transm. Problems, ~
        e-mail: blank@iitp.ru}
        \thanks{This research has been partially supported
                by Russian Foundation for Basic Research
                and program ONIT.}
       }
\date{May 20, 2011} 
\maketitle

\begin{abstract}%
We study deterministic discrete time exclusion type spatially
heterogeneous particle processes in continuum. A typical example
of this sort is a traffic flow model with obstacles: traffic
lights, speed bumps, spatially varying local velocities etc.
Ergodic averages of particle velocities are obtained and their
connections to other statistical quantities, in particular to
particle and obstacles densities (the so called Fundamental
Diagram) is analyzed rigorously. The main technical tool is a
``dynamical'' coupling construction applied in a nonstandard
fashion: instead of proving the existence of the successful
coupling (which even might not hold) we use its presence/absence
as an important diagnostic tool. %
\end{abstract}%

\section{Introduction}\label{s:intro}

In 1970 Frank Spitzer introduced the (now classical) simple
exclusion process as a Markov chain that describes
nearest-neighbor random walks of a collection of particles on the
one-dimensional
infinite\footnote{or finite with periodic boundary conditions} %
integer lattice. Particles interact through the hard core
exclusion rule, which means that at most one particle is allowed
at each site. This seemingly very particular process appears
naturally in a very broad list of scientific fields starting from
various models of traffic flows \cite{NS,GG,ERS,Bl-erg,Bl-hys},
molecular motors and protein synthesis in biology, surface growth
or percolation processes in physics (see \cite{Pe,BFS} for a
review), and up to the analysis of Young diagrams in
Representation Theory \cite{CMOS}.

From the point of view of the order of particle interactions there
are two principally different types of exclusion processes: with
synchronous and asynchronous updating rules. In the latter case at
each moment of time a.s. at most one particle may move and hence
only a single interaction may take place. This is the main model
considered in the mathematical literature (see e.g. \cite{Lig} for
a general account and \cite{An,EFM,IK} for recent results), and
indeed, the assumption about the asynchronous updating is quite
natural in the continuous time setting. The synchronous updating
means that {\em all} particles are trying to move simultaneously
and hence an arbitrary large (and even infinite) number of
interactions may occur at the same time. This makes the analysis
of the synchronous updating case much more difficult, but this is
what happens in the discrete time case.\footnote{if one do not
   consider some ``artificial'' updating rules like a sequential
   or random updating.} %
This case is much less studied, but still there are a few results
describing ergodic properties of such processes
\cite{Bl-erg,Bl-hys,ERS,GG,NS}.

Recently in \cite{B09} we have introduced and studied the
synchronous updating version of the simplest exclusion process
(TASEP) in continuum.\footnote{Some other continuous space
   generalizations were considered e.g. in \cite{KG,ZS}.}  %
Now our aim is to extend these results to the case when the media
is heterogeneous, i.e. contains obstacles. The idea is as follows.
Consider a (countable) collection of particles performing (random)
walks on a real line with hard core type interactions (we assume
that the particles cannot outran each other). Assume also that on
the real line there is a (countable) number of obstacles. To
overcome an obstacle a particle needs to spend some additional
time on it. Thus we have two types of interactions: between
particles and with static obstacles. The principal novelty here is
that the presence of obstacles leads to a spontaneous creation of
``traffic jams'' (particles clusters) near obstacles. Indeed, it
is easy to construct initial particle configurations having no
clusters of particles but such that these clusters will be created
in front of obstacles under dynamics. This feature was completely
absent in the case without obstacles: a cluster may only disappear
but never be created. We start the formalization with the simplest
version of the model and later in
Sections~\ref{s:varying},\ref{s:gen} consider what happens in a
more complicated setting.

By a {\em configuration} $x:=\{x_i\}_{i\in\IZ}$ we mean an ordered
(i.e. $x_i\le x_{i+1}~\forall i\in\IZ$) bi-infinite sequence of
real numbers $x_i\in\IR$. The union $X$ of all such sequences
plays the role of the phase space of the system under
consideration. Consider also a special configuration $z\in X$ such
that $z_i\le z_{i+1}~\forall i\in\IZ$. We refer to elements of $x$
as to positions of particles and to elements of $z$ as to
positions of obstacles.

\def\particle#1{
     \put(0,0){\circle{3}} \bline(0,0)(0,-1)(25) \put(-3,-32){$x_{#1}$}
     \put(0,5){\vector(1,0){40}} \put(20,9){$\xi_{#1}$}}
\Bfig(180,60)
      {\footnotesize{
       \thicklines
       \bline(0,10)(1,0)(180)
       \put(30,10){\particle{i}} \put(150,10){\particle{i+1}}
       \thinlines
       \bline(30,-10)(1,0)(90) \put(65,-8){$\t\G_i$}
       \bline(120,-16)(0,1)(40) \put(115,-20){$z_j$}
       \bline(30,30)(1,0)(120)  \put(80,33){$\G_i$}
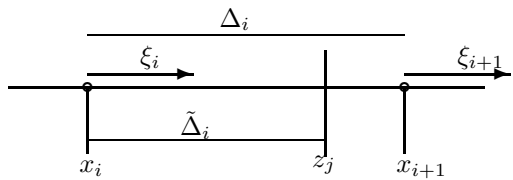
      }      } {TASEP in heterogeneous continuum. \label{f:tasep-c}}

Let $\G_i=\G_i(x):=x_{i+1}-x_i$ and let %
$$ \t\G_i=\t\G_i(x,z):=\min\left(x_{i+1}-x_i,
                ~~\min_j(z_j:~z_j>x_i) - x_i\right) $$
stands for the minimum between the distances from the point
$x_i^t$ to the point $x_{i+1}^t$ and to the next obstacle (see
Fig.~\ref{f:tasep-c}). We refer to $\G_i,\t\G_i$ as {\em gaps} and
{\em modified gaps} corresponding to the $i$-th particle in $x$.

For a given real ${v}>0$ (to which we refer as a maximal local
velocity) and the configuration $z$ we define a map $\map:X\to X$
as follows: %
\beq{e:map}{(\map x)_i := x_i + \min\left(\t\G_i(x,z),{v}\right)
                        \qquad\forall i\in\IZ.} %
The dynamical system $(\map,X)$ describes the collective motion of
particles discussed above. If $\t\G_i(x,z)<v$ we say that the
$i$-th particle is {\em blocked} (meaning that its motion is
blocked by the $i+1$-th particle or an obstacle) at time $t$ and
{\em free} otherwise. By a {\em cluster} of particles in a
configuration $x\in X$ we mean consecutive particles with gaps
$\G_i(x)<v$ having no obstacles between them.

Associating the number of iterations of the map $\map$ with time
we often use the notation $\map^tx\equiv x^t:=\{x_i^t\}$. In this
terms the quantity $\xi_i^t:=\min\left(\t\G_i(x^t,z),{v}\right)$
plays the role of the {\em local velocity} of the $i$-th particle
in the configuration $x\equiv x^0$ at time $t\ge0$ and thus the
dynamics can be rewritten as %
\beq{e:dynamics}{x_i^{t+1}:=x_i^{t}+\xi_i^t .} %

A more general setting including varying/random velocities and
waiting times will be considered in
Sections~\ref{s:varying},\ref{s:gen}.

It is of interest that even without obstacles (i.e. if
$z=\emptyset$) the behavior of the deterministic dynamical system
$(\map,X)$ is far from being trivial. In \cite{B09} it was shown
that this system is chaotic in the sense that its topological
entropy is positive (and even infinite).

Our main results are concerned with the so called Fundamental
Diagram describing the connection between average particle
velocities $V(x)$ and particle/obstacle densities $\rho(x),
\rho(z)$ (see Section~\ref{s:metric} for definitions) and
technically the analysis is based on a (somewhat unusual)
``dynamical'' coupling construction (see Section~\ref{s:couping}
and also \cite{B09}).

For a given $v>0$ and a configuration of obstacles $z$ denote by
$\t{z}$ the {\em extended} configuration of obstacles obtained by
inserting between each pair of entries $z_i,z_{i+1}$ new
$\intp{(z_{i+1}-z_i)/v}$ `virtual' obstacles at distances
$v$ between them starting from the point $z_i$. Here
$\intp{u}$ stands for the integer part of the number $u$.

\begin{theorem}\label{t:fund-diagram} (Fundamental Diagram)
Let $\rho(x), \rho(\t{z})$ be well defined. Then %
\beq{e:vel}{ V(x)=\min\{1/\rho(\t{z}),~1/\rho(x)\} .}%
Therefore the phase space is divided into two parts: gaseous %
$\{(x,\t{z}):~\rho(x)\le\rho(\t{z})\}$ (consisted of
configurations of eventually non-interacting particles) and liquid
$\{(x,\t{z}):~\rho(x)>\rho(\t{z})\}$ (where clusters of particles
are present at any time).
\end{theorem}

\Bfig(320,100)
      {\footnotesize{
       \put(0,0){\bpic(0,0)
       \thicklines
       \bline(0,0)(1,0)(130) \bline(0,0)(0,1)(90)
       \put(132,-5){$\rho(x)$}  \put(-17,85){$\rho(\t{z})$}
       \thinlines
       \bline(0,0)(3,2)(130)
       \put(30,55){{\bf gas}} \put(60,15){{\bf liquid}} }
       \put(190,0){\bpic(0,0)
       \thicklines
       \bline(0,0)(1,0)(140) \bline(0,0)(0,1)(90)
       \bline(0,73)(1,0)(60)
       \bezier{200}(60,73)(75,5)(135,2)
       \thinlines
       \bezier{30}(60,73)(60,36)(60,0)
       \put(145,-5){$\rho(x)$}  \put(-8,85){$V(x)$}
       \put(-15,70){$\frac1{\rho(\t{z})}$}
       \put(85,30){$1/\rho(x)$} \put(-8,-5){$0$}
       \put(55,-12){$\rho(\t{z})$}
      }      }      }
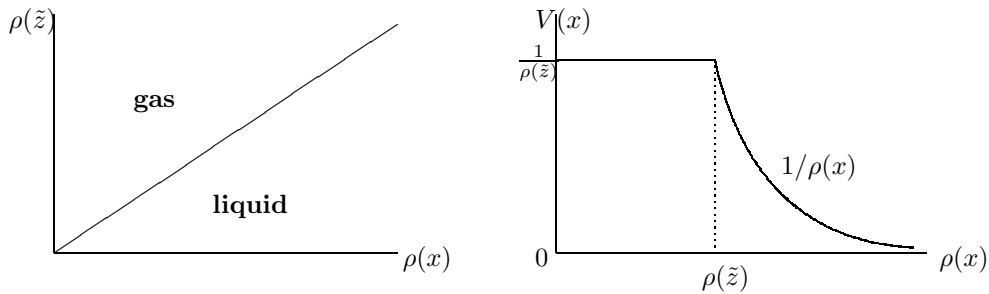
{(left) Phase Diagram, (right) Fundamental Diagram.
\label{f:tasep}}

These results are illustrated in Fig~\ref{f:tasep}. It is worth
note that a naive argument tells that the average velocity of a
particle in an infinite system is essentially the minimum of the
average inter-particle distance and the average inter-obstacle
distance. Theorem~\ref{t:fund-diagram} shows that this is
absolutely not the case: instead of the average inter-obstacle
distance one needs to take into account the average distance
between the elements of the extended configuration $\t{z}$. As we
shall see the latter is very different from the former, e.g. we
always have $\rho(\t{z})\ge1/v$ independently on $\rho(z)$ (which
even might be not well defined). Actually the mere fact that the
correspondence between the average velocities and particle/obstacle
densities is one-to-one comes as a surprise especially without any
regularity assumptions for the positions of obstacles.

The main steps of the proof are as follows. First we show that for
given configurations of particles and obstacles upper/lower
average particle velocities are the same for all particles. Now to
compare average velocities in different particle configurations
with the same configuration of obstacles we develop a special
dynamical coupling construction. Applying this construction we
prove that that upper/lower average particle velocities are the
same provided the same particle densities. Thus to calculate the
dependence of average velocities on densities it is enough to
construct for each density a single configuration having this
density for which we are able to calculate its velocity
explicitly. To this end we consider an auxiliary zero-range
lattice process whose trajectories are invariant under the
dynamics of our original system without obstacles. To calculate
the average velocity in the zero-range lattice process one uses
corresponding results obtained in \cite{B09}.

Despite that various couplings are widely used in the analysis of
Interacting Particle Systems (IPS) (see e.g. \cite{Lig}),
applications of our approach is very different from usual. In
particular, we do not prove the existence of the so called
successful coupling (which even might not hold) but instead use
its presence/absence as an important diagnostic tool. Remark also
that typically one uses the coupling argument to prove the
uniqueness of the invariant measure and to derive later other
results from this fact. In our case there might be a very large
number of ergodic invariant measures or no invariant measures at
all (e.g. the trivial example of a single particle performing a
skewed random walk). This example reminds about another important
statistical quantity -- average particles velocity. The dynamical
coupling will be used directly to find connections between the
average particle velocities and other statistical features of the
systems under consideration, in particular with the corresponding
particle densities.

It is worth note that all approaches used to study discrete time
lattice versions of IPS are heavily based on the combinatorial
structure of particle configurations. This structure has no
counterparts in the continuum setting under consideration. In
particular the particle -- vacancy symmetry is no longer
applicable in our case. This explains the need to develop a
fundamentally new techniques for the analysis of IPS in continuum.
The presence of obstacles also prevents the direct application of
the scheme developed in \cite{B09} for the spatially homogeneous
case.

\?{This techniques cannot be applied directly in the lattice case.
Nevertheless, the embedding of lattice systems to the continuum
setting allows to obtain (indirectly) new results for the lattice
systems as well.}

\section{Basic properties}\label{s:metric}

Here we shall study questions related to particle densities and
velocities.

By the {\em density}
$$ \rho(x,I):= \frac{\#\{i\in\IZ:~~x_i\in I\}}{|I|} $$
of a configuration $x\in X$ in a bounded segment $I=[a,b]\in\IR$
we mean the number of particles from $x$ whose centers $x_i$
belong to $I$ divided by the Lebesgue measure $|I|>0$ of the
segment $I$. If for any sequence of {\em nested} bounded segments
$\{I_n\}$ with $|I_n|\toas{n\to\infty}\infty$ the limit
$$\rho(x):=\lim\limits_{n\to\infty}\rho(x,I_n)$$
is well defined we call it the {\em density} of the configuration
$x\in X$. Otherwise one considers upper and lower (with respect to
all possible collections of nested intervals $I_n$) particle
densities $\rho_\pm(x)$. %

\begin{remark}\label{r:density} If $\exists \rho(x)<\infty$ then
$|x_n-x_m|/|n-m|\toas{|n-m|\to\infty}1/\rho(x)$.
\end{remark}

\begin{lemma}\label{l:density-preservation}
The upper/lower densities $\rho_\pm(x^t)$ are preserved by
dynamics, i.e. $\rho_\pm(x^t)=\rho_\pm(x^{t+1})~~\forall t\ge0$.
\end{lemma}%
\proof For a given segment $I\in\IR$ the number of particles from
the configuration $x^t\in X$ which can leave it during the next
time step cannot exceed 1 and the number of particles which can
enter this segment also cannot exceed 1. Thus the total change of
the number of particles in $I$ cannot exceed 1, because if a
particle leaves the segment through one of its ends no other
particle can enter through this end. Therefore $$|\rho(x^t,I) -
\rho(x^{t+1},I)|\cdot|I|\le1$$ which implies the claim. \qed%

By the (average) {\em velocity} of the $i$-th particle in the
configuration $x\in X$ at time $t>0$ we mean
$$V(x,i,t):=\frac1t\sum\limits_{s=0}^{t-1}\xi_i^s %
         \equiv(x_i^{t}-x_i^0)/t.$$ %
If the limit %
$$V(x,i):=\lim\limits_{t\to\infty}V(x,i,t)$$ is well defined we call
it the (average) {\em velocity} of the $i$-th particle. Otherwise
one considers upper and lower particle velocities $V_\pm(x,i)$.

\?{Denote by $\t{z}$ the {\em extended} configuration of obstacles
$z$ defined as follows. First we set $\t{z}:=z$. Then $\forall
i\in\IZ$ if $\t{z}_i+{v}<\t{z}_{i+1}$ we add the point
$\t{z}_i+{v}$ to $\t{z}$ and re-enumerate $\t{z}$. This procedure
continues until $\t{z}_{i+1}-\t{z}_i\le{v}~~\forall i\in\IZ$.
Eventually we get the configuration satisfying the property
$0<\t{z}_{i+1}-\t{z}+i\le v~\forall i\in\IZ$.}

\begin{lemma}\label{l:velocity-preservation}
Let $x\in X$ and let $\rho(\t{z})<\infty$ be well defined. Then
$|V(x,j,t) - V(x,i,t)|\toas{t\to\infty}0$ a.s. $\forall
i,j\in\IZ$.
\end{lemma}%
\proof It is enough to prove this result for $j=i+1$. Consider the
difference between (average) velocities of consecutive particles%
\bea{ V(x,i+1,t) - V(x,i,t)
 \a= \frac{x_{i+1}^t-x_{i+1}^0}t - \frac{x_i^t-x_i^0}t \\ %
 \a= \frac{x_{i+1}^t-x_i^t}t - \frac{x_{i+1}^0-x_i^0}t \\ %
 \a= \G_i^t/t - \G_i^0/t .} %
The last term vanishes as $t\to\infty$ and it is enough to show
that the same happens with $\G_i^t/t$. Here and in the sequel we
use the notation $\G_i^t\equiv\G_i(x^t)$.

Normally in the deterministic setting gaps between particles
$\G_i^t$ are uniformly bounded (see e.g.
\cite{Bl-erg,Bl-hys,B09}). Surprisingly in the present setting
this is not the case and we are only able to show that $\G_i^t$ may
grow not faster than $o(t)$, which fortunately is enough for our
aims.

To prove this estimate we introduce a new configuration having
only two particles $\2{x}:=\{\2{x}_1,\2{x}_2\}$ which have the
same initial positions as the $i$-th and the $(i+1)$-th particles
of $x$, i.e. $\2{x}_1=x_i^0, \2{x}_2=x_{i+1}^0$. Denoting
$\2{x}^t:=\map^t\2{x}$ and $\2{\G}^t:=\2{x}_2^t-\2{x}_1^t$ we
observe that $\2{\G}^t\ge \G_i^t~~\forall t\ge0$.

Consider the movement of the leading particle in the process
$\t{x}$. There is the first moment of time $t_2:=t_2(\t{x}_2)$
when this particle encounters an obstacle from $z$. Denote by
$\t{z}_k$ the position of this obstacle in the extended
configuration $\t{z}$. By the definition of the dynamics
$\map^t\t{z}_k=\t{z}_{k+t}~~\forall t\ge0$.

The movement of the trailing particle in the process $\t{x}$ is a
bit more complicated since additionally to obstacles it may be
stopped by the leading particle. Anyway there is the first moment
of time $t_1:=t_1(\t{x}_1,\t{x}_2)>t_2$ when this particle
encounters the obstacle located at $\t{z}_k$. For each $t>t_1$
these two particles move synchronously along the ``chain''
$\t{z}$. Thus the growth of $\2{\G}^{t}$ is completely determined
by the concentration of entries in $\t{z}$ (i.e. by
$\rho(\t{z})$). Due to the assumption about the existence of the
density of $\t{z}$ these fluctuations cannot exceed $o(t)$ and
hence $\G_i^t/t\le\2{\G}^t/t\toas{t\to\infty}0$. \qed

\begin{corollary}\label{c:vel-equiv}
Under the assumptions of Lemma~\ref{l:velocity-preservation} the
upper and lower particle velocities $V_\pm(x,i)$ do not depend on
$i$.
\end{corollary}

\begin{remark} (a) If the density $\rho(\t{z})$ is not well defined,
the distance $\G_i^t$ may grow linearly with time and the limit
average velocities might differ for different particles. (b) The
existence of $\rho(z)<\infty$ does not imply the existence of
$\rho(\t{z})<\infty$ but fir ``typical'' $z$ this implication
takes place (see Section~\ref{s:gen}).
\end{remark}

Let us estimate upper/lower densities of the extended
configuration $\t{z}$ of obstacles for a given configuration $z$
with $\rho(z)>0$ and a given local velocity $v>0$.

\begin{lemma}\label{l:extended} %
$\max\{1/v,\rho(z)\} \le \rho_-(\t{z})
                     \le \rho_+(\t{z}) \le \rho(z) + 1/v$.
\end{lemma}
\proof Observe that for each $i\in\IZ$ the spatial segment
$[z_i,z_{i+1}]$ contains exactly %
$\lfloor(z_{i+1}-z_i)/v\rfloor\in[(z_{i+1}-z_i)/v-1,(z_{i+1}-z_i)/v]$
elements. Therefore $\forall n\in\IZ_+$ %
\bea{ \rho(\t{z},[z_0,z_n])
  \a= \left( n+\sum_{i=0}^{n-1}\lfloor(z_{i+1}-z_i)/v\rfloor \right)
   /(z_n-z_0) \\
  \a= \rho(z,[z_0,z_n])
  + \left( \sum_{i=0}^{n-1}\lfloor(z_{i+1}-z_i)/v\rfloor \right)
   /(z_n-z_0) \\
  \a\le \rho(z,[z_0,z_n]) + 1/v .} %
Similarly %
$$ \rho(\t{z},[z_0,z_n]) \ge \rho(z,[z_0,z_n]) - n/(z_n-z_0) + 1/v
                         = 1/v .$$
Passing to the limit as $n\to\infty$ and using that $\t{z}$
contains all elements of $z$ we get the result. \qed

\section{Dynamical coupling through particle's overtaking} \label{s:couping}

Consider two independent particle processes $x^t,\2x^t$. %
By $i(x^t)$ we denote the $i$-th particle of the process $x^t$
(i.e. $x_i^t$ is the location of the particle $i(x^t)$ at time
$t$).

We say that the particle $i(x^t)$ {\em overtakes} at time $t>0$
the particle $j(\2x^t)$ if $x_i^{t-1}<\2x_j^{t-1}$ and
$x_i^t\ge\2x_j^t$, and denote this event as
$i(x^t)\hookrightarrow{j(\2x^t)}$.

\begin{lemma}\label{l:over} If \/ $i(x^t)\hookrightarrow{j(\2x^t)}$
then \\ %
(a) $\not\exists n\in\IZ: ~ n(x^t)\hookrightarrow{j(\2x^t)}$; \\ %
(b) $\not\exists m\in\IZ: ~ j(\2x^t)\hookrightarrow{m(x^t)}$. \\
Additionally \?{the event} $i(x^t)\hookrightarrow j(\2x^t)$ and
$(j-1)(\2x^t)\hookrightarrow i(x^t)$ might happen only if
$x_i^t=\2x_{j-1}^t$.
\end{lemma}

Note that there is no contradiction between two parts of this
Lemma since the former part concerns a particle being overtaken
while the latter concerns a overtaking particle.

\proof By the construction of the process under consideration and
the definition of the overtaking we have: %
\beq{e:over}{x_{i-1}^t \le x_i^{t-1} < \2x_j^{t-1} \le \2x_j^{t}
  \le x_i^t \le x_{i+1}^{t-1} .} %
Now (a) follows from the observation that by (\ref{e:over}) neither
particles preceding $i(x^t)$ nor particles succeeding it may
overtake the particle $j(\2x^t)$ at time $t$.

If (b) would hold then %
$\2x_j^{t-1} < x_m^{t-1} \le x_m^t \le \2x_j^t$. Thus by
(\ref{e:over}) we get $\2x_j^{t-1} \le x_m^{t-1} \le \2x_j^t$,
which might be possible only if $\2x_j^{t-1}=\2x_j^t$. The latter
contradicts to the assumption about the overtaking.

To show that the event discussed in the additional part may take
place, consider a pair of configurations satisfying for some
$i,j,t\in\IZ$ the following inequalities: %
$$ \2x_{j-1}^{t-1} < x_i^{t-1} < \2x_j^{t-1}
                   = \2x_{j+1}^{t-1} = x_{i+1}^{t-1}, \qquad
  \max\left(\t\Delta_i(x^{t-1},z),
            \t\Delta_{j-1}(\2x^{t-1},z)\right) \le v .$$
Then $x_i^t=\2x_{j-1}^t=\2x_j^t$ which implies that %
$i(x^t)\hookrightarrow j(\2x^t)$ and $(j-1)(\2x^t)\hookrightarrow
i(x^t)$. Assume now that on the contrary $x_i^t\ne\2x_{j-1}^t$. If
$x_i^t>\2x_{j-1}^t$ then $\2x_{j-1}^t\le\2x_{j}^t<x_i^t$ which
contradicts to $(j-1)(\2x^t)\hookrightarrow i(x^t)$. Similarly
$x_i^t<\2x_{j-1}^t$ implies $x_i^t<\2x_{j-1}^t\le\2x_{j}^t$ which
contradicts to $i(x^t)\hookrightarrow j(\2x^t)$. \qed

Let us introduce the {\em dynamical coupling} of the processes
$x^t,\2x^t$ which consists in a consequent pairing of particles of
opposite processes. The pairing in our deterministic setting is a
pure formal action which does not change the dynamics. The idea is
that if a particle overtakes some particles from the opposite
process it becomes paired with one of them. The construction
starts at time $t=0$ and initially all particles are assumed to be
unpaired.

Denote by
$$ J_i(x^t):=\{j\in\IZ: ~ i(x^t)\hookrightarrow{j(\2x^t)}\} $$
the set of particles overtaken simultaneously by the particle
$i(x^t)$ at time $t$, by
$$ i_\bullet(x^t):=\min\left\{\infty,
~\min\{j\in J_i(x^t): ~ j(\2x^{t-1}) ~{\rm is ~paired} \}
\right\}$$ the paired particle with the minimal index among them,
and finally
$$ i_\circ(x^t):=
   \function{i_\bullet(x^t)-1  &\mbox{if } i_\bullet(x^t)-1\in J_i(x^t) \\
             \infty &\mbox{otherwise} }. $$%

\?{if $i_\bullet(x^t)-1\in J_i^t$ we set
$i_\circ(x^t):=i_\bullet(x^t)-1$ and $i_\circ(x^t):=\infty$
otherwise.}

In words, $i_\bullet(x^t)$ stands for the paired particle in
$J_i^t$ having the minimal index, and $i_\circ(x^t)$ stands for
the unpaired particle in $J_i(x^t)$ having the maximal index among
those preceding $i_\bullet(x^t)$. A typical connection between
those indices and positions of the $i$-th particle at times $t-1$
and $t$ is shown in Fig~\ref{f:index}.

\Bfig(180,60)
      {\footnotesize{
       \thinlines 
       \bline(0,0)(1,0)(180) \bline(0,30)(1,0)(180)
       \put(5,30){\circle{4}} \put(135,30){\circle{4}}
       \put(155,30){\circle*{4}} \put(175,30){\circle*{4}}
       \put(20,0){\circle{4}} \put(40,0){\circle{4}}
       \put(100,0){\circle*{4}} \put(115,0){\circle*{4}}
       \put(3,35){$x_i^{t-1}$}  \put(133,35){$x_i^{t}$}
       \put(38,6){$\2x_{i_\circ}^{t}$}
       \put(98,6){$\2x_{i_\bullet}^{t}$}
       \put(-10,28){$x$} \put(-10,-2){$\2x$}
      }      }
      {Connections between indices and particle's positions. \label{f:index}}

To simplify the description of the dynamical coupling we shall use
a diagrammatic representation for coupled configurations, where
paired particles are denoted by black circles and unpaired ones by
open circles, and use the upper line of the diagram for the
$x$-particles (i.e. particles from the $x$-process) and the lower
line for the $\2x$-particles. In this representation a typical
pairing event looks as follows \\ %
\beq{e:diag}{
_{\circ~~~~~~~\circ\circ}^{~~~~\circ\circ}\Longrightarrow~
 _{~~\circ~~~\circ~~~~~\circ}^{\circ~~~~~~\circ}~\longrightarrow~ %
 _{~~\bullet~~~\star~~~~~\circ}^{\bullet~~~~~~\star} .} %
Here $\Longrightarrow$ corresponds to the dynamics and
$\longrightarrow$ to the pairing, while
~$_\bullet^\bullet,~_\star^\star$~ stand for two different mutual
pairs of particles created after the particle overtaking under
dynamics.

Now we are ready to define the pairing rigorously. We proceed
first with all $x$-particles, overtaking some $\2x$-particles at
time $t$, and then with all $\2x$-particles, overtaking some
$x$-particles at this time.

Let the overtaking takes place for $i$-th $x$-particle at time
$t$, then if %
\begin{itemize}
\item[(a)] $i(x^{t-1})$ is paired and $i_\circ(x^t)<\infty$ we
           re-pair these particles: %
$$_{\bullet~~\circ\circ\star\circ}^{~\bullet~~~\star}~\Longrightarrow~ %
_{\bullet~\circ\circ\star~~~\circ}^{~~~~\bullet~~\star}~\longrightarrow~ %
_{\circ~\circ\bullet\star~~~\circ}^{~~~~\bullet~~\star} $$

\item[(b)] $i(x^{t-1})$ is unpaired and %
 \begin{itemize}
 \item[(b')] $i_\bullet(x^t)<\infty$ and $x_i^t>\2x_{i_\bullet(x^t)}$
         we re-pair these particles: %
         $$_{\circ~~\bullet\star}^{~~\circ~~\bullet\star}~\Longrightarrow~ %
          _{\circ\bullet~~\star}^{~~\circ\bullet~~\star}~\longrightarrow~ %
          _{\circ\bullet~~\star}^{~~\bullet\circ~~\star}$$
 \item[(b'')] else if $i_\circ(x^t)<\infty$ this unpaired particle forms
          a new pair with $i(x^{t})$: 
          $$_{~~\circ\circ\star}^{\circ~~~~\star}~\Longrightarrow~ %
           _{\circ\circ~~\star}^{~~\circ~~\star}~\longrightarrow~ %
           _{\circ\bullet~~\star}^{~~\bullet~~\star}$$
 \end{itemize}
\end{itemize}

The pairing rules when the overtaking takes place for $i$-th
$\2x$-particle are exactly the same except for the exchange of $x$
by $\2x$ and vice versa.

The complexity of these rules reflects that first a particle may
overtake simultaneously several particles from another process,
and second an arbitrary number of particles may share the same
spatial position. In particular, in the following event we have: %
$$_{~~~\circ~~\bullet\star}^{\circ~~~~~\bullet}~\Longrightarrow~~ %
 _{\circ~~\bullet~~~\star}^{~~~\circ~~~\bullet}~\longrightarrow~
 _{\bullet~~\bullet~~~\star}^{~~~\bullet~~~\bullet} $$
rather than %
$~\longrightarrow~
 _{\circ~~\bullet~~~\star}^{~~~\bullet~~~\circ}$. %
\?{(Here we use the same notation as in diagram (\ref{e:diag}).)}
Indeed, the trailing $x$-particle overtakes both two trailing
$\2x$-particles and since the re-pairing rules in (b) demand
strict inequalities only the pairing of trailing unpaired
particles takes place.

By Lemma~\ref{l:over} each overtaking event may be treated
separately. Therefore we apply the pairing rules first for all
$x$-particles, overtaking some $\2x$-particles, and then for all
$\2x$-particles, overtaking some $x$-particles.

According to the definition, particles from the same pair move
synchronously until either the movement of one of them is blocked
or one of the members of the pair is swapped with an unpaired
particle from the same process. It is convenient to think about
the coupled process as a ``gas'' of single (unpaired) particles
and ``dumbbells'' (pairs). A previously paired particle may
inherit the role of the unpaired one from one of its neighbors. In
order to keep track of positions of unpaired particles we shall
refer to them as $x$- and $\2x$-{\em defects} depending on the
process they belong.

We say that a pair of configurations $(x,\2x)$ is {\em proper} if
for each two mutually paired particles the open segment between
them cannot exceed ${v}$, does not contain neither a defect nor an
obstacle (i.e. the situations $_\bullet^{~\circ~\bullet}$,
$_\bullet|^\bullet$ cannot happen), and there are no crossing
mutually paired pairs ($_{\star~\bullet}^{~\bullet~\star}$).

\begin{lemma}\label{l:pairing-proper} If the pair of configurations
$(x^{t-1},\2x^{t-1})$ is proper for some $t>1$, then the pair
$(x^t,\2x^t)$ is proper as well.
\end{lemma}
\proof The situation $_\bullet^{~\circ~\bullet}$ may happen if
either a defect overtakes the trailing particle in a pair, or if a
pair is born around this defect. Both such possibilities
contradict to the definition of pairing.

Assume that at time $t$ there is a pair of mutually paired
particles $i(x^t), j(\2x^t)$ such that the open segment between
them contains an obstacle: $x_i^t<z_k<\2x_j^t$. Assume also that
this pair was present at time $t-1$ as well. According to the
construction of the dynamics this may happen only if
$x_i^{t-1}<z_k-{v}$ and $z_k=\2x_j^{t-1}$, which contradicts to
the assumption that the pair $(x^{t-1},\2x^{t-1})$ is proper:
namely the distance between members of the same pair exceeds
${v}$. On the other hand if the particles $i(x^t), j(\2x^t)$ were
not paired at time $t-1$ and the pair is just created at time $t$
then one of the particles should overtake another at this time.
Hence the distance between these particles cannot exceed ${v}$.

It remains to show that the last property still holds for ``old''
pairs. In order to enlarge the distance between the existing
mutually paired particles one of them should be blocked by another
particle or by an obstacle. On the other hand, the obstacle cannot
belong to the open segment between pair members and hence it may
block only the leading particle in the pair, which may only
decrease the distance.

The blocking particle might be paired or unpaired. The former case
implies that the `companion' of the blocking particle is at
distance at most $v$ and hence it will block the enlargement of
the distance by more than $v$. In the letter case the non-blocked
paired particle overtakes the unpaired one and hence the pair
under consideration will be re-paired.

The observation that in the moment of the creation of a pair the
distance between its members cannot exceed $v$ finishes this part
of the proof.

The analysis of the absence of crossing mutually paired pairs is
completely similar to the absence of defects between elements of a
pair and therefore we skip this point. \qed

Denote by $\rho_u(x,I)$ the density of the $x$-defects belonging
to a finite segment $I$, and by $\rho_u(x):=\rho_u(x,\IR)$ the
upper limit of $\rho_u(x,I_n)$ taken over {\em all} possible
collections of nested finite segments $I_n$ whose lengths go to
infinity.

We say that a coupling of two Markov particle processes
$x^t,\2x^t$ is {\em nearly successful} if the upper density of the
$x$-defects $\rho_u(x)$ vanishes with time a.s.
This definition differs significantly from the conventional
definition of the successful coupling (see e.g. \cite{Lig}), which
basically means that the coupled processes converge to each other
in finite time.

\begin{lemma}\label{l:velocity-coupled} Let $x,\2x\in X$ with
$\rho(x)=\rho(\2x)>0$, and let there exist a nearly successful
dynamical coupling $(x^t,\2x^t)$. Then
$$|V(x,0,t)-V(\2x,0,t)|\toas{t\to\infty}0.$$
\end{lemma}%
\proof Consider an integer valued function $n_t$ which is equal to
the index of the $\2x$-particle paired at time $t>0$ with the
$0$-th $x$-particle. If the $0$-th $x$-particle is not paired at
time $t$ we set %
$n_t:=\function{n_{t-1} &\mbox{if } t>0 \\
                0       &\mbox{if } t=0}$. %

To estimate the growth rate of $|n_t|$ at large $t$ observe that
$n_t$ changes its value only at those moments of time when the
$0$-th $x$-particle meets a $\2x$-defect. By the assumption about
the nearly successful coupling at time $t\gg1$ the average
distance between the defects at time $t$ is of order
$1/\rho_u(\2x^{t})$ while the amount of time needed for two
particles separated by the distance $L$ to meet cannot be smaller
than $L/(2{v})$. Therefore the frequency of interactions of the
$0$-th $x$-particle with
$\2x$-defects may be estimated from above by the quantity of order %
$\rho_u(\2x^t)\toas{t\to\infty}0$, which implies %
$n_t/t\toas{t\to\infty}0$.

Now we are ready to prove the main claim. %
\bea{ |V(x,0,t) - V(\2x,0,t)| \a= |(x_0^t-x_0^0) -
(\2x_0^t-\2x_0^0)|/t \\
\a\le |x_0^t-\2x_0^t|/t + |x_0^0-\2x_0^0|/t \\%
\a\le |x_0^t-\2x_{n_t}^t|/t +
        \frac{|n_t|}t~|\2x_{n_t}^t-\2x_0^t|/|n_t|
                            + |x_0^0-\2x_0^0|/t .} %
The 2nd addend is a product of two terms $|n_t|/t$ and
$|\2x_{n_t}^t-\2x_0^t|/|n_t|$. As we have shown, the 1st of them
vanishes with time. If $|n_t|$ is uniformly bounded, then the 2nd
term is obviously uniformly bounded on $t$. Otherwise, for large
$|n_t|$ by Remark~\ref{r:density} and the density preservation the
2nd term is of order $1/\rho(\2x)$, which proves its uniform
boundedness as well. Thus the 2nd addend goes to 0 as
$t\to\infty$. Noting finally that the 1st and the last addend also
vanishes with time at rate $1/t$ we are getting the result. \qed

\begin{lemma}\label{l:s-coupling}
Let $\rho(x)=\rho(\2x)$ and let in the coupled process $\forall
i,j~~\exists$ a (random) moment of time $t_{ij}<\infty$ such that
$x_i^t>\2x_j^t$ for each $t\ge t_{ij}$. Then the coupling is
nearly successful.
\end{lemma}
\proof By the assumption each $x$-particle will overtake
eventually each $\2x$-particle located originally to the right
from its own position and thus will form a pair with it or with
one of its neighbors (if they are so close that were overtaken
simultaneously). Thus the creation of pairs is unavoidable. To
show that the upper density of defects cannot remain positive,
consider how the defects move under our assumptions. Assume that
at time $t\ge0$ the $i$-th $x$-particle is paired with the $j$-th
$\2x$-particle. Then by Lemma~\ref{l:pairing-proper} in order to
overtake at time $s>t$ the $j$-th $\2x$-particle significantly (by
a distance larger than $v$) the $i$-th $x$-particle necessarily
needs to break the pairing with the $j$-th $\2x$-particle. Thus by
the definition of the dynamical coupling either a $x$-defect
overtakes the
$j$-th $\2x$-particle: %
~$_{~\bullet}^{\circ~~~\bullet}~\longrightarrow~
  _\bullet^{~\circ~\bullet}~\longrightarrow~
  _\bullet^{~\bullet~\circ}$,
or the $i$-th $x$-particle overtakes a $\2x$-defect:
~$_{\bullet~~~\circ}^{~~~\bullet}\longrightarrow~
  _{\bullet~\circ}^{~~~\bullet}\longrightarrow~
  _{\circ~\bullet}^{~~~\bullet}$. %
(Otherwise this pair will not be broken.) Therefore during this
process the $x$-defects move to the right while the $\2x$-defects
move to the left. Hence they inevitably meet each other and
``annihilate''. The assumption about the equality of particle
densities implies the result. \qed

\section{Auxiliary lattice zero-range process}\label{s:lattice}

Consider now a lattice process which we shall need in the sequel.
This process is defined on an integer lattice $\IZ$ occupied by a
bi-infinite configuration of particles $y$ ordered with respect to
their positions $y_i$, i.e. $y_i\le y_{i+1}~~\forall i$. The
dynamics is defined as follows. For each $i\in\IZ$ consider all
particles occupied the site $i$ and choosing the one with the
largest index (say $k_i$) among them we move this particle by one
position to the right, i.e. $y_{k_i}:=y_{k_i}+1$. This is a
deterministic version of the so called zero-range process on $\IZ$
with parallel updating rules illustrated in
Fig~\ref{f:zero-range}.

\Bfig(180,60)
      {\footnotesize{
       \thinlines 
       \bline(0,0)(1,0)(180)
       \put(20,2){\circle{4}} \put(20,7){\circle{4}}
       \put(20,12){\circle{4}} \put(20,17){\circle{4}}
       \put(20,22){\circle{4}}
       \put(20,2){\vector(4,1){70}} \put(0,20){\vector(4,1){20}}
       \put(90,2){\circle{4}} \put(90,7){\circle{4}}
       \put(90,12){\circle{4}} \put(90,17){\circle{4}}
       \put(90,2){\vector(2,1){70}}
       \put(160,2){\circle{4}} \put(160,7){\circle{4}}
       \put(160,12){\circle{4}} \put(160,17){\circle{4}}
       \put(160,22){\circle{4}} \put(160,27){\circle{4}}
       \put(160,32){\circle{4}} 
       \put(160,2){\vector(2,1){20}}
       \put(10,-10){$i-1$} \put(90,-10){$i$} \put(150,-10){$i+1$}
      }      }
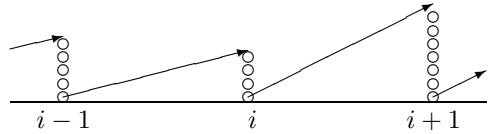
      {Zero-range process. \label{f:zero-range}}

Setting $v=1, ~z=\emptyset$ and assuming that $x_i\in \IZ~~\forall
i\in\IZ$ we observe that in this case (without obstacles) the
trajectory $\map^t x$ coincides with the trajectory of the
zero-range process. Therefore according to \cite{B09} the average
particle velocity for the zero-range process is equal to %
\beq{e:v-hom}{ V(y)=\min\{1,~1/\rho(y)\} .} %

Assume now that the lattice under consideration is not uniform
(e.g. $\IZ$) but the distance between the $i$-th and $(i+1)$-th
site is equal to a nonnegative number $\ell_i$ for each $i\in\IZ$.
Considering the zero-range process on this heterogeneous lattice,
but assuming that the site's occupation is the same in both cases
we are able to calculate the corresponding statistical quantities.

Denote by $\2y$ a configuration on a heterogeneous lattice $\2\IZ$
in which distances between the $i$-th site and the $(i+1)$-th one
are given by the sequence of numbers $\{\ell_i\}$. Let $y$ be a
configuration on $\IZ$ such that the $i$-th particle of the
configuration $\2y$ is located on the site whose index coincides
with the index of the site occupied by the $i$-th particle of the
configuration $y$.

\begin{lemma}\label{l:heter-l}
The deterministic zero-range processes on $\IZ$ and $\2\IZ$
starting with configurations $y,\2y$ described above for each
$t\ge0$ preserves the connection between the particle
configurations. Therefore the particle velocities in
these processes satisfy the relation %
\beq{e:vel-heter}{ \2V(\2y)=V(y)/\rho(\2\IZ) ,}  %
where $\rho(\2\IZ)$ is assumed to be positive and is defined as
the density of the particle configuration having exactly one
particle at each site of~ $\2\IZ$.
\end{lemma}
\proof The first claim follows from the definition of the
zero-range process. Denote $L(t):=y_0^t-y_0^0~~\forall{t}$. Then %
$$ (\2y_0^t-\2y_0^0)/t
 = \frac1t \sum_{j=0}^{y_0^t-y_0^0-1}\ell_{j+y_0^0} %
 = \frac{L(t)}t~ \frac1{L(t)}\sum_{j=0}^{L(t)-1}\ell_{j+y_0^0} \\ %
 \toas{t\to\infty} V(y)/\rho(\2z). $$
\qed

\section{Proof of Theorem~\ref{t:fund-diagram}}

\begin{lemma} Let $x,\2x\in X$ $\rho(x)=\rho(\2x)$ and let $V(x)$
be well defined for given $v,z$. Then $V(\2x)$ is also well
defined and $V(\2x)=V(x)$.
\end{lemma}%
\negvskip%
\proof Let $x,\2x\in X_\rho:=\{z\in X:~~\rho(z)=\rho\}$ be two
admissible configurations of the same particle density. If one
assumes that the dynamical coupling procedure leads to the nearly
successful coupling of particles in these configurations then by
Lemma~\ref{l:pairing-proper} the assumptions of
Lemma~\ref{l:velocity-coupled} are satisfied and hence
$|V(x,0,t)-V(\2x,0,t)|\toas{t\to\infty}0$ which by
Lemma~\ref{l:velocity-preservation} implies the claim.
In general the assumption about the nearly successful coupling may
not hold,\footnote{Consider e.g. the setting with $1/\rho>5v$ and
     the configurations $x_i:=i/\rho$ and $\2x_i:=i/\rho+2v$, and
     assume that there are no obstacles. Then $\rho(x)=\rho(\2x)=\rho$,
     $V(x)=V(\2x)=v$ but no pair will be created.} %
however as we demonstrate below the pairing construction is still
useful.

Define random variables
$$W_{ij}^t:=x_i^t-\2x_j^t,~i,j\in\IZ,~t\in\IZ_0.$$ Then %
$$V(x,i,t)-V(\2x,j,t)=W_{ij}^t/t - W_{ij}^0/t.$$ Since by
Lemma~\ref{l:velocity-preservation} the differences between
average velocities of different particles belonging to the same
configuration vanish with time it is enough to consider only the
case $i=j=0$. For $W_{00}^t$ there might be three possibilities
which we study separately:

\begin{itemize}
\item[(a)] $\lim\limits_{t\to\infty} W_{00}^t/t=0$. Then %
$|V(x,0,t)-V(\2x,0,t)|\le|W_{00}^t|/t + |W_{00}^0|/t
 \toas{t\to\infty}0$, %
which by Corollary~\ref{c:vel-equiv} implies that the sets of
limit points of the average velocities coincide.

\item[(b)] $\limsup\limits_{t\to\infty}W_{00}^t/t>0$. Then
$\forall i\in\IZ$ the $i$-th particle of the $x$-process will
overtake eventually each particle of the $\2x$-process located at
time $t=0$ to the right from the point $x_i^0$. This together with
the assumption of the equality of particle densities allows to
apply Lemma~\ref{l:s-coupling} according to which the coupling is
nearly successful. On the other hand, by
Lemma~\ref{l:pairing-proper} the distance between mutually paired
particles cannot exceed $v$. Therefore by
Lemma~\ref{l:velocity-coupled} we have ~~%
$|V(x,0,t)-V(\2x,0,t)|\toas{t\to\infty}0$, %
which contradicts to the assumption (b).

\item[(c)] $\limsup\limits_{t\to\infty} W_{00}^t/t<0$. Changing
the roles of the processes $x^t,\2x^t$ one reduces this case to
the case (b).
\end{itemize}

Thus only the case (a) may take place. \qed

To apply this result to prove Theorem~\ref{t:fund-diagram} one
needs to construct for each particle density $\alpha$ a
configuration having this density for which we are able to
calculate its velocity explicitly. There are two possibilities to
realize this idea. In both cases we use the auxiliary lattice
zero-range process $y^t$ constructed in Section~\ref{s:lattice} on
the heterogeneous lattice $\2\IZ$ whose $i$-th site coincides with
the location of the $i$-th element in the extended configuration
of obstacles $\t{z}$.

The first construction is as follows. Choose an arbitrary initial
configuration of particles $y$ of a given density for the
zero-range process. By Lemma~\ref{l:heter-l} and the relation
(\ref{e:v-hom}) we obtain the relation~(\ref{e:vel}). Observing
that the trajectory of the zero-range process $y^t$ coincides with
the trajectory of our original process $\map^ty$ we get the
desired result.

An alternative way to derive the relation (\ref{e:vel}) is to
construct a specific initial configuration of the zero-range
process, for which we can calculate the corresponding average
velocity directly. The key observation here is that for $x:=\t{z}$
obviously we have $V(x)=1/\rho(\t{z})$.

For a given $x\in X$ for which $\exists~0<\rho(x)<\infty$ and
$\alpha\in\IR_+$ we define the configuration $\alpha x\in X$ in
three steps depending on arithmetic properties of $\alpha$: \par%
(a)~$kx$ for $\alpha=k\in\IZ_+$; \par%
(b)~$\frac{k}n x$ for $\alpha=\frac{k}n$, $k=mn+k'$
    and $m,k'\in\IZ_0$, $k'<n\in\IZ_+$; \par%
(c)~$\frac{k_n}n x\toas{n\to\infty}\alpha x$ if
$\frac{k_n}n\toas{n\to\infty}\alpha$. %

Here by $kx$ we mean a configuration in which each particle in $x$
is replaced by $k$ particles. To get $y:=\frac{k}n x$ we divide
$x$ into blocks having $n$ consecutive particles each and
construct $y:=\{y_i\}$ as follows. For each $i$ at location $x_i$
we set exactly $m$ particles at positions $y_j=x_i$. Besides at
positions of each of the first $k'$ particles in each of the
blocks we add an additional particle and enumerate the particles
according to their positions.

If $\alpha\ne k/n$ one considers a sequence of its rational
approximations: $k_n/n\toas{n\to\infty}\alpha$. For each $n$
according to the rule (b) we construct a configuration
$x^{(\alpha_n)}:=\frac{k}nx$. Observe that for $\ell>1$ the
configuration $x^{(\alpha_{n\ell})}$ differs from $x^{(\alpha_n)}$
only by the presence of some additional particles in each of the
blocks of length $n\ell$. Using the existence of $\rho(x)$ one
proves that the limit configuration exists, and we denote the
latter by $\alpha x$.


Now for any $\alpha\in\IR_+$ we construct $x:=\alpha\t{z}$ and
calculate $V(\alpha\t{z})$ as follows.
We start with the case $\alpha=1$, i.e. $x=\t{z}$. In this case
obviously $\map x$ coincides with the left shift of the
configuration $x$. Therefore %
$(\t{z}_t-\t{z}_0)\cdot\rho(\t{z},[\t{z}_0,\t{z}_t])\equiv t
~~\forall t\ge0$
and thus %
$$ V(x,0,t)=(\t{z}_t-\t{z}_0)/t = 1/\rho(\t{z},[\t{z}_0,\t{z}_t])
   \toas{t\to\infty}1/\rho(\t{z}) .$$

If $\alpha=\frac{k}n$ the configuration $x:=\alpha\t{z}$ is
``almost'' spatially periodic: it consists of blocks of the
configuration $\t{z}$ having equal number of $n$ particles. The
spatial periodicity immediately implies the existence of the
average velocity $V(x)$. Observe now that for any rational
$\alpha<1$ the particles in the configuration $\alpha\t{z}$ on
average move exactly as in the original configuration $\t{z}$
(i.e. $V(\alpha\t{z})\equiv V(\t{z})$). If $\alpha=\frac{k}n>1$
then the direct calculation gives $V(x)=\frac{n}kV(\t{z})$.
Passing to the
limit as $\frac{k_n}n\toas{n\to\infty}\alpha$ we get eventually %
$V(\alpha\t{z}):=\function{
       \frac1{\rho(\t{z})} &\mbox{if ~~ } \alpha\le1 \\
       \frac1{\alpha\rho(\t{z})} &\mbox{otherwise }}$. %
Here one uses that $V(\frac{k}n\t{z})\ge
V(\frac{k}{n\ell}\t{z})~~\forall\ell\in\IZ_+$.

\section{Varying waiting times and local velocities}\label{s:varying}

In this Section we consider a more general model which includes
both varying waiting times and local velocities. Namely we assume
that the waiting time at the $j$-th obstacle is equal to
$\tau_j\in\{1,2,\dots,\tau\}$, and the local velocity for
particles in the spatial segment $(z_j,z_{j+1})$ is equal to
$v_j\in(0,v]$. It is useful to think about this model as a road
divided into parts separated by obstacles (e.g. traffic lights)
and having different quality of the pavement, which we describe by
varying from part to part local velocities.

If $\tau_j>1$ it might be possible that a new particle is coming
to an obstacle when some preceding particles are still waiting
there. We assume that for the given particle the waiting time
$\tau_j$ at the $j$-th obstacle starts only when the succeeding
particle leaves it. The movement of a particle waiting at a
certain obstacle in the next moment of time depends on the time
which this and preceding particles already spent there. Therefore
this model is non Markovian. To cure this pathology using methods
developed in \cite{Bl-hys} (for a very different situation) one
adds a new variable {\em type} (for each particle) whose value is
equal to the amount of time the particle will wait at an obstacle
if it is located on the obstacle and zero otherwise.

Strictly speaking in order to use the notion of the coupling this
is important (since it is defined for a pair of Markov processes).
Nevertheless in the deterministic setting which we consider in
this paper this issue is not crucial (as we shall see) and we
proceed without this Markovian extension.

Denote by $\b{\tau},\b{v}$ the collections of waiting times and
local velocities corresponding to corresponding obstacles. If
$\tau_j\equiv\tau$ and $v_j\equiv{v}$ we recover the previous
setting. Surprisingly the analysis of this general setting is very
similar to the previous one. Therefore we shall discuss only the
points when the analysis differs.

First we refine the configuration of obstacles $z$ exchanging the
original $i$-th obstacle by $\tau_i+1$ consequent obstacles
located at the same position $z_i$ but having 0 waiting times,
i.e. leaving one these new obstacles the particle immediately
moves to the next of them. The resulted ordered collection of
obstacles we again call $z$. To take care about the change of
indices we change also the collection of local velocities
inserting $\tau_i$ unit velocities (before the original element
$v_i$) corresponding to new obstacles and re-enumerating them.

Obviously for any configuration of particles $x$ the movements of
its elements in the case of the original configurations of
obstacles and velocities and refined ones are exactly the same.

To construct the {\em extended configuration} of obstacles $\t{z}$
we proceed as follows: between the obstacles $z_i$ and $z_{i+1}$
we insert $\lfloor(z_{i+1}-z_i)/v_i\rfloor$ virtual obstacles at
distance $v_i$ starting from the point $z_i$ for each $i\in\IZ$.
If $v_i\equiv v$ this construction boils down to the one described
in Section~\ref{s:intro}.

It is straightforward to check that all constructions and results
obtained in Sections~2-4 remain valid in this more general setting
except only one point. In the definition of the proper pair of
configurations one changes that the open segment between two
mutually paired particles belonging to the interval
$[z_i,z_{i+1}]$ for some $i\in\IZ$ cannot exceed $v_i$ (instead of
$v$ in the original definition).

Therefore using the same arguments as in the proof of
Theorem~\ref{t:fund-diagram} we get its generalization.

\begin{theorem}\label{t:fund-diagram-} (Fundamental Diagram)
Let $z,~\tau_j\in\{1,2,\dots,\tau\},~v_j\in(0,v]$ be given and
let $\rho(x), \rho(\t{z})$ be well defined and positive. Then %
\beq{e:vel1}{ V(x)=\min\{1/\rho(\t{z}),~1/\rho(x)\} .}%
\end{theorem}

\section{Discussion and generalizations}\label{s:gen}

\def\bnum#1{\n{\bf{#1}}}

\bnum1. The density of a configuration in the way how it was
defined in Section~\ref{s:metric} depends sensitively on the
statistics of both left and right tails of the configuration. A
close look shows that in fact if all particles move in the same
direction, say right, one needs only the information about the
corresponding (right) tail, which allows to expand significantly
the set of configurations having densities and for which our
results can be applied.

For a configuration $x\in X$ by a {\em one-sided particle density}
we mean the limit %
\beq{e:one-side}
    {\hat\rho(x):=\lim_{\ell\to\infty}\rho(x,[0,\ell]).}%
The upper an lower one-sided densities correspond to the upper and
lower limits.

\begin{theorem}\label{t:one-side} All previous results formulated in
terms of ``two-sided'' densities remain valid if one replaces the
usual particle density $\rho$ to the one-sided density $\hat\rho$.
\end{theorem}
\proof The key observation here is that the movement of a given
particle in a configuration $x^t\in X$ depends only on particles
with larger indices. Therefore if one changes positions of all
particles with negative indices the particles with positive
indices will still have the same average velocity. On the other
hand, by Lemma~\ref{l:velocity-preservation} the average velocity
does not depend on the particle index. This allows to apply the
following trick.

For each configuration $x\in X$ of density $\rho(x)$ we associate
a new configuration $\hat{x}\in X$ defined by the relation:
$$ \hat{x}_i:=\function{x_i^t &\mbox{if } i\ge0 \\
                        x_0+i/\rho(x) &\mbox{otherwise }.} $$
Then obviously $\hat\rho(x)=\rho(\hat{x})=\rho(x)$.

Therefore for all purposes related to the average velocities all
results valid for the configuration $\hat{x}$ remain valid for $x$
as well. \qed

\bnum2. A close look to the proof of Theorem~\ref{t:fund-diagram}
may lead to the hypothesis that for a given configuration of
obstacles $z$ ``typical'' trajectories eventually become supported
only by the locations belonging to the extended configuration
$\t{z}$. Indeed, the calculation of the average velocity is given
exactly for the configurations of this type. Let us show that this
is not the case even in the simplest setting when the
configuration of obstacles has density one and is supported by
integer points (i.e. a single obstacle at each point of $\IZ$).
Let $\rho(x)\ge2$ and let the initial configuration has at least 2
particles at each half-integer point $\frac12~\IZ$. Then for each
$v\ge1/2$ and $t\ge1$ the configuration $x^t$ is supported by the
lattice $\IZ \cup \frac12~\IZ$.

\bigskip

\bnum3. Question about the existence of $\rho(\t{z})$. Let
$\rho(z)>0$ be well defined. For each given $v>0$ we say that the
configuration $z$ is {\em regular} if the density of the
corresponding extended configuration $\t{z}$ is well defined. In
the topology induced by the uniform metric in the space of
sequences there exists an open set containing irregular
configurations $z$. Nevertheless assuming a reasonable model for
the creation of the configuration of obstacles one can show that a
``typical'' configuration of obstacles is regular.

Indeed, assume that the sequence $z$ is a realization of an
ergodic stochastic process with stationary increments. Then
Birkhoff Ergodic Theorem implies that Cesaro means for the
sequence of fractional parts $\{(z_{i+1}-z_i)/v\}$ converge to the
limit and thus our claim follows.

If $\rho_+(\t{z})\ne\rho_-(\t{z})$ then one may consider
upper/lower particle velocities and the corresponding statement
for the Fundamental Diagram may be rewritten as follows.

\begin{theorem}\label{t:fund-diagram-h}
Let $v,\rho(z)>0$ be given. Then %
$$ V_\pm(x)=\min\{1/\rho_\mp(\t{z}),~1/\rho_\mp(x)\}  .$$
\end{theorem}

The proof of this result follows basically the same arguments as
in the case of Theorem~\ref{t:fund-diagram} but some additional
technical estimates related to partial limits are necessary.

\bigskip

\bnum4. Assume now that we consider our original setting (i.e.
$\tau_j\equiv0~\forall j\in\IZ$) and the local velocity of the
$i$-th particle does not depends on space (as in
Section~\ref{s:varying}) but depends on time according to a given
collection of local velocities $\{v_i^t\}$, i.e. $v_i^t$ stands
for the local velocity of the $i$-th particle at time $t$. Thus
using the notation introduced in Section~\ref{s:intro} we get
$\xi_i^t:=\min\left(\t\Delta_i(x,z),v_i^t\right)$ and
$x_i^{t+1}:=x_i^{t}+\xi_i^t$.

\begin{theorem}\label{t:dif-vel} For any given configuration of
obstacles $z$ such that $z_k\toas{k\to\infty}\infty$ and two
one-particle configurations $x:=\{x_0\}\ne\2x:=\{\2x_0\}$ there
exists a sequence of local velocities $\{v_0^t\}$ such that
$|x_0^t-\2x_0^t|\ge Ct$ for all $t\in\IZ_+$ and some $C>0$. Hence
average particle velocities cannot coincide.
\end{theorem}%

The proof is based on the observation that one can choose
$\{v_0^t\}$ such that each time $t$ when the $x$-particle meets an
obstacle it makes a ``full'' step equal to $\{v_0^t\}$ while the
$\2x$-particle which meets with a different obstacle makes only a
half step.

This result shows that a direct generalization of our pure
deterministic setting to the random case is not available.

\newpage%

\end{document}